\def\a{\alpha}

\def\b{\beta}

\def\brel{\buildrel}
\def\bsk{\bigskip}

\def\CC{{\bf C}}

\def\cC{{\cal C}}

\def\cE{{\cal E}}
\def\cF{{\cal F}}

\def\cJ{{\cal J}}

\def\cl{\colon}
\def\cL{{\cal L}}
\def\cM{{\cal M}}

\def\cMh{\widehat{\cal M}}
\def\cMt{\widetilde{\cal M}}
\def\cN{{\cal N}}
\def\cO{{\cal O}}
\def\cod{{\rm cod}}

\def\cQ{{\cal Q}}

\def\cS{{\cal S}}
\def\cU{{\cal U}}
\def\cV{{\cal V}}

\def\d{\delta}
\def\D{\Delta}

\def\deg{{\rm deg}}

\def\det{{\rm det}}

\def\e{\epsilon}

\def\es{\emptyset}

\def\g{\gamma}

\def\G{\Gamma}

\def\hb{\hbox}

\def\Hom{{\rm Hom}}
\def\hra{\hookrightarrow}

\def\i{\iota}

\def\Im{{\rm Im}}

\def\l{\lambda}

\def\L{\Lambda}
\def\lra{\longrightarrow}

\def\msk{\medskip}
\def\n{\noindent}

\def\om{\omega}
\def\omt{\widetilde{\omega}}
\def\op{\oplus}
\def\ot{\otimes}
\def\ov{\overline}
\def\O{\Omega}
\def\Oh{\wh{\Omega}}

\def\pf{\noindent{\bf Proof.}\hskip 2mm}
\def\PGL{{\rm PGL}}
\def\phit{\widetilde{\phi}}
\def\Pic{{\rm Pic}} 

\def\pih{\widehat{\pi}}
\def\pit{\widetilde{\pi}}
 
\def\PP{{\bf P}}

\def\qed{\hfill{\bf q.e.d.}} 
\def\QQ{{\bf Q}}

\def\rk{{\rm rk}}

\def\s{\sigma}

\def\Si{\Sigma}
\def\Sih{\wh{\Sigma}}
\def\Sit{\wt{\Sigma}}
\def\sm{\setminus}

\def\ss{\subset}

\def\t{\theta}

\def\tB{\widetilde{B}}

\def\tm{\times}

\def\tV{\widetilde{V}}

\def\vf{\varphi}

\def\wh{\widehat}
\def\wt{\widetilde}

\def\ZZ{{\bf Z}}
\centerline{\bf DESINGULARIZED MODULI SPACES OF SHEAVES ON A $K3$, II.}  
\bsk
\centerline{\bf May 20 1998}
\msk
\centerline{\bf Kieran G.~O'Grady}
\bsk
\bsk
\n
{\bf 0. Introduction.}
\bsk
This is a sequel to~[O4]. In that paper we constructed a symplectic
desingularization of $\cM_4$, the moduli space of rank-two semistable
torsion-free  sheaves $F$ on a $K3$ surface $X$ with $c_1(F)=0$ and
$c_2(F)=4$. In this paper we will prove that our desingularization, denoted
$\cMt_4$, is a new (ten-dimensional) irreducible symplectic variety.
Explicitely, we will show that $\cMt_4$ is one-connected and that
$h^{2,0}(\cMt_4)=1$; this means that $\cMt_4$ is an irreducible symplectic
variety. Furthermore we will prove that $b_2(\cMt_4)\ge 24$. Since all
known ten-dimensional irreducible symplectic varieties have $b_2=7$ or
$b_2=23$~[B,H], our results show that all deformations of $\cMt_4$ are new
irreducible symplectic varieties. The paper is organized as follows. In the first (and
longest) section we prove that $\cMt_4$ is connected and that $h^{2,0}(\cMt_4)=1$: 
we adapt to our situation Jun Li's strategy for determining the stable $b_1,b_2$
of moduli spaces of sheaves on a surface~[Li2]. In the second section we prove $\cMt_4$ is
simply-connected: first we prove that if $\cO_X(1)$ has degree $2$ (we can assume this by
a deformation argument), $\cMt_4$ is birational to a certain Jacobian fibration over
$|\cO_X(2)|$, then we show by a monodromoy argument  that the Jacobian fibration is
simply-connected. In the third section we exhibit a $24$-dimensional subspace of
$H^2(\cMt_4;\QQ)$. 
\bsk
\bsk
\n
{\bf 1.  $\cMt_4$ is irreducible and $h^{2,0}(\cMt_4)=1$.} 
\bsk
\n
The proof goes as follows. First, by
well-known arguments involving deformations of polarized $K3$'s~[GH,O1 \S 2] we can
assume the polarization has genus two, i.e.~$H^2=2$; this will be assumed throughout the
paper. We will also assume that $X$ is a very general genus two $K3$, i.e.~that
$Pic(X)=\ZZ[H]$; thus $H$ is $c$-generic for any $c$~[O4,(0.2)], 
so that by~[O4, \S 2] there is a symplectic desingularization $\pit\cl\cMt_4\to\cM_4$. We
consider  the morphism $\phi\cl\cM_4\to\PP^N$  associated to a
high power of the determinant line-bundle~[LP1,Li1], and we let 
$\phit\cl\cMt_4\to\PP^N$ be the  composition with the desingularization
map $\pit$. We will show that $\phit$ is semi-small, hence we can apply
the Lefschetz Hyperplane Section (LHS) Theorem~[GM] as if $\phit$ were an
embedding. Let  $C\in|H|$ be a smooth curve: following Jun Li we
choose a particular linear subspace $\L\ss\PP^N$ of codimension at most
$4$, such that   
$$\phit^{-1}\L=\tV_C\cup\Sit_C\cup\tB_C,\leqno(1.0.1)$$   
where $\tV_C,\Sit_C,\tB_C$ are closed subsets such that:
\msk
\item{(1)}
$\tV_C\sm(\Sit_C\cup\tB_C)$ consists of the points $x$ such that $\pit(x)$ parametrizes a
locally-free sheaf whose restriction to $C$ is not semistable.
\item{(2)}
$\Sit_C=\pit^{-1}\Si_C$, where 
$$\Si_C:=\{[I_Z\op I_W]|\ \hb{$Z$ or $W$ intersects $C$}\}.$$
\item{(3)}
$x\in\tB_C\sm\Sit_C$ if and only if  $\pit(x)$ parametrizes a stable sheaf which is not
locally-free at some point of $C$.    
\msk
\n
This is where the hypothesis that $C$ has genus two is used: since the
moduli space of rank-two semistable vector-bundles on $C$ with trivial
determinant has dimension three, the intersection of four generic
theta-divisors is empty, and this gives Decomposition~(1.0.1) for a $\L$ of
codimension at most four. By the LHS Theorem, the map induced by
inclusion      
$$H^q(\cMt_4;\ZZ)\to H^q(\phit^{-1}\L;\ZZ)$$ 
is an isomorphism for $q\le 5$. Thus we are reduced to analyzing $\tV_C$,
$\Sit_C$, $\tB_C$, and how they intersect.
The component $\tV_C$ will be described in terms of
$\cM(H,3)$, the moduli space of semistable rank-two torsion-free sheaves $E$ on $X$ with
$c_1(E)=H$, $c_2(E)=3$. In fact, if $x\in\tV_C$ is generic, and $F$ is the sheaf
parametrized by $\pit(x)$, let $G$ be the elementary modification of $F$ associated to the
desemistabilizing sequence for $F|_C$, and let $E:=G(1)$; then $E$ is a stable rank-two
vector-bundle with $c_1(E)=H$, $c_2(E)=3$, hence $[E]\in\cM(H,3)$. Conversely $F$ can
be reconstructed from $E$ and the choice of a certain rank-one subsheaf of $E|_C$.
Equivalently, $\tV_C$ is birational to a moduli space of parabolic sheaves 
$F_0=F\supset F_1=G\supset F_2=F(-C)$. Building on our knowledge of
$\cM(H,3)$, we will show that $\tV_C$ is birational to
$\Pic^1(C)\tm\cM(H,3)\tm\PP^1$.  On the other hand $\Sit_C$, $\tB_C$ 
are birational to $\PP^1$-fibrations over $C\tm X\tm X^{(2)}$ and
$C\tm X^{(3)}$, respectively. We will also give birational
descriptions of components of the intersection between $\Sit_C$ and
$\tV_C$, $\tB_C$. At this point we will be ready to prove the main results
of the section. First, since       
$$\tV_C\cap\Sit_C\not=\es\not=\Sit_C\cap\tB_C,$$ 
we get $b_0(\cMt_4)=1$. Next, if $V'_C$,
$\Si'_C$, $B'_C$ are desingularizations of $\tV_C$, $\Sit_C$, $\tB_C$ respectively, the
pull-back map  
$$\a\cl\G(\O^2_{\cMt_4})\to\G(\O^2_{V'_C})\op\G(\O^2_{\Si'_C})\op\G(\O^2_{B'_C})$$
is injective by the LHS Theorem. Of course the components of a point in $\Im\a$
must satisfy some compatibily conditions: these conditions will imply that $\Im\a$ is
one-dimensional, hence $h^{2,0}(\cMt_4)=1$.     

\proclaim Remark.
{\rm The procedure outlined above should allow us to determine the Betti
numbers of $\cMt_4$ up to $b_5$ included. However, even the determination of $b_2$ (which
should be equal to $24$) requires a much more detailed analysis of
$\tV_C\cup\Sit_C\cup\tB_C$; the calculation of $h^{2,0}$ is simpler because $h^{2,0}$ is a
birational invariant.}

\bsk
\n
{\bf 1.1. The determinant map $\phi$.}
\bsk
\n
Le Potier and Jun Li~[LP1,Li1] have constructed a determinant line bundle  $\cL$ on $\cM_4$
(in general on any moduli space of semistable torsion-free sheaves on a surface). They
proved that if $m\gg 0$ then $\cL^{\ot m}$ is base-point free~[Li1, Thm.3] and that the
{\it determinant map}  
$$\phi_m\cl\cM_4\to\PP\left(H^0(\cM_4,\cL^{\ot m})^{*}\right)$$
can be identified with the map to the Uhlenbeck compactification~[Li1 Thm.4,FM].
From now on we fix  a very large $m$ and we set $\phi=\phi_m$. We will give an explicit
description of $\phi$.  Let $\cM_4^{st},\cM_4^{lf}\ss\cM_4$ be the (open)  subsets
parametrizing stable and locally-free sheaves, respectively.  By Lemma~(1.1.5) of~[O4],
$$\cM_4\sm\cM_4^{st}=\Si:=\{[I_Z\op I_W]\ |
\ \ \hb{$Z,W$ zero-dimensional subschemes of $X$ with $\ell(Z)=\ell(W)=2$}\},$$ 
hence $\cM_4^{lf}\ss\cM_4^{st}$. In particular~[Li1, Thm.4] the restriction of $\phi$ to 
$\cM^{lf}$ is an isomorphism onto its image. Let the {\it boundary} of $\cM_4$
be $B:=\cM_4\sm\cM_4^{lf}$, and set $B^{st}:=B\cap\cM_4^{st}$. 

\proclaim (1.1.1) Proposition.
Keeping notation as above, let $[F]\in B^{st}$. Then  $F^{**}\cong\cO_X^{(2)}$.

\pf
Since $F$ is not locally-free, $c_2(F^{**})<c_2(F)=4$. Thus Hirzebruch-Riemann-Roch gives
$\chi(F^{**})>0$. By Serre duality $h^2(F^{**})=h^0(F^{*})=h^0(F^{**})$,
hence  $h^0(F^{**})>0$. A non-zero section of $F^{**}$ must have isolated zeroes by
slope-semistability, hence we have an exact sequence  
$$0\to\cO_X\brel\a\over\to F^{**}\to I_Z\to 0,$$
where $\ell(Z)=c_2(F^{**})$. If $\ell(Z)=0$ then $F^{**}\cong\cO_X^{(2)}$, and we are
done. By Cayley-Bacharach $\ell(Z)\not=1$. Now assume $\ell(Z)\ge 2$: then we can find a
subsheaf $I_W\ss\cO_X$, with $\ell(W)=(4-\ell(Z))$, such that $\a(I_W)\ss F$, and this
contradicts stability of $F$. 
\qed
\msk
In order to describe $\phi$ on the boundary we introduce the {\it singularity cycle} of a point
$[F]\in\cM_4$: we assume that the representative $F$ has been chosen so that it
is isomorphic to the direct sum of the successive quotients of its Harder-Narasimhan
filtration: this means that if $F$ is strictly semistable, then $F\cong I_Z\op I_W$. With this
hypothesis we set 
$$sing[F]:=\sum\limits_{p\in X}\ell_F(p),$$
where $\ell_F(p)$ is the length at $p$ of the Artinian sheaf $F^{**}/F$. It follows from 
Lemma~(1.1.5) of~[O4] and Proposition~(1.1.1)  that  if $[F]\in B$, then $sing
[F]\in X^{(4)}$. By~[Li1, Thm.4] the restriction of $\phi$ to $B$ is simply the singularity map
$$\matrix{ B & \brel\phi\over\lra & X^{(4)} \cr
                   [F] & \mapsto & sing[F] \cr}.$$
As is easily checked the image $\phi(B)$ is equal to the whole $X^{(4)}$, hence  
the Uhlenbeck stratification~[FM] becomes
$$\Im\phi=\cM_4^{lf}\coprod X^{(4)}.\leqno(1.1.2)$$
\bsk
\n
{\bf 1.2. The map $\phit$.}
\bsk
\n
Let $\phit:=\phi\circ\pit$, where $\pit\cl\cMt_4\to\cM_4$ is the 
desingularization of~[O4, \S 2]. Let
$$\G_d:=\{x\in\Im\phit|\ \dim\phit^{-1}(x)=d\}.$$

\proclaim (1.2.1) Proposition.
Keep notation as above. The map $\phit$ is semi-small, i.e.
$$\cod(\G_d,\Im\phit)\ge 2d$$
for all $d$.  

Roughly speaking the reason is that the non-degenerate two-form $\omt_4$ is in a
"stratified" sense the pull-back of a two-form on $\Im\phit$ (strictly speaking this is
nonsense because  $\Im\phit$ is singular). More precisely: there is a stratification 
$$X^{(4)}=\coprod\limits_{p} X^{(4)}_p\leqno(1.2.2)$$ 
indicized by partitions of $4$. If $p=(p_1,p_2,p_3,p_4)$ is a partition,
i.e.~$\sum\limits_{i=1}^4ip_i=4$, the stratum $X^{(4)}_p$ parametrizes cycles
$$\g=\sum\limits_{j=1}\limits^{p_1}q_j^1+2\sum\limits_{j=1}\limits^{p_2}q_j^2+\cdots,$$
where the $q_j^i$ are pairwise distinct. Letting $X_{\d}^{(s)}:=\left(X^{(s)}\sm
diagonals\right)$, we have an open inclusion
$$X^{(4)}_p\ss X_{\d}^{(p_1)}\tm X_{\d}^{(p_2)}\tm\cdots,$$
in particular $X^{(4)}_p$ is smooth. Let
$\om^{(s)}\in\G(\O^2_{X_{\d}^{(s)}})$ be the holomorphic two-form
obtained from $\om$ by symmetrization, and for a partition $p$,  let
$\om_p^{(4)}$ be the two-form on $X^{(4)}_p$ given by  
$$\om_p^{(4)}:=\sum\limits_{i=1}\limits^4 i \pi_i^*\om^{(p_i)},$$
where $\pi_i\cl X^{(4)}_p\to X_{\d}^{(p_i)}$ is the $i$-th projection. Stratification~(1.2.2)
gives a stratification of $\Im\phit=\Im\phi$, in which $\cM_4^{lf}$ is the open stratum, and
the $X_p^{(4)}$ are the remaining strata. The stratified pull-back formula is the following:
$$\omt_4|_{\phit^{-1}X_p^{(4)}}=-{1\over 4\pi^2}\phit^{*}\om_p^{(4)}.\leqno(1.2.3)$$
It is an immediate consequence of~[O2,(2-9)]. 
\msk
\n
{\bf Proof of Proposition~(1.2.1).}
\hskip 2mm
Equivalently,  we must show that
$$\cod(\phit^{-1}\G_d,\cMt_4)\ge d.\leqno(1.2.4)$$
Let $V\ss \phit^{-1}\G_d$ be an irreducible component. If
$V\cap\cM_4^{lf}\not=\es$, then $d=0$ and there is nothing to prove. So assume
$V\ss B=\phit^{-1}X^{(4)}$. Then $V\cap\phit^{-1}X_p^{(4)}$  is dense in $V$ for some
partition $p$. Let $y\in V$ be such that $V$ is smooth  at $y$ and the restriction $\phit\cl
V\to \phit(V)$ is submersive at $y$; it follows from~(1.2.3) that 
$$\omt_4\left(T_yV,T_y\phit^{-1}(\phit(y))\right)\equiv 0.$$
Thus $\omt_4$ induces a map 
$$T_y\cMt_4/T_yV\to T_y\phit^{-1}(\phit(y))^{*},$$
which is surjective because $\omt_4$ is non-degenerate. Hence
$$\cod(V,\cMt_4)=\dim\left(T_y\cMt_4/T_yV\right)\ge\dim T_y\phit^{-1}(\phit(y))=d.$$
This proves~(1.2.4).
\qed

\proclaim Corollary.
Let $p=(p_1,\ldots,p_4)$ be a partition of $4$. If $\g\in X_p^{(4)}$, then
$$\dim\phit^{-1}(\g)\le 5-\sum\limits_{i=1}\limits^4 p_i.\leqno(1.2.5)$$

\pf
The map 
$$\phit^{-1} X_p^{(4)}\brel\phit\over\lra X_p^{(4)}$$
is a locally trivial fibration in the analytic topology, hence the fibers
have a constant dimension, say $d$. Thus $X_p^{(4)}\ss\G_d$, so that
$$\cod(\G_d,\Im\phit)\le\cod(X_p^{(4)},\Im\phit).$$
By Proposition~(1.2.1) we get
$$2\dim\phit^{-1}(\g)=2d\le\cod(X_p^{(4)},\Im\phit).$$
Writing out the codimension on the right, one gets Inequality~(1.2.5).
\qed 
\msk
Let $h^0(\cM_4,\cL^{\ot m})=N+1$, so that $\phit\cl\cMt_4\to\PP^N$. 

\proclaim (1.2.6) Corollary.
Let $\L\ss\PP^N$ be a linear subspace of codimension at most $c$. The map
$$H^q(\cMt_4;\ZZ)\to H^q(\phit^{-1}\L;\ZZ)$$
induced by inclusion is an isomorphism for all $q\le (9-c)$. 

\pf
An immediate consequence of Proposition~(1.2.1) and the (generalized) LHS Theorem~[GM,
p.150]. 
\qed
\bsk
\n
{\bf 1.3. Choice of $\L$.}
\bsk
\n
We fix once and for all a smooth curve $C\in |H|$. Let $\t$ be line bundle on $C$ of degree one
(i.e.~half the canonical degree); thus for $[F]\in\cM_4$
$$\chi(F|_C\ot\t)=0.\leqno(1.3.1)$$
There is a canonical section~[FM,Li1] $\s_{\t}$ of the determinant line-bundle $\cL$ such
that 
$$supp(\s_{\t})=\{[F]\in\cM_4|\ h^0(F|_C\ot\t)>0\}.$$
By~(1.3.1) the right-hand side has codimension at most one in every component of $\cM_4$,
but a priori it might contain a whole component. In any case there exists a linear subspace
$\L_{\t}\ss\PP^N$ of codimension at most one such that
$\phi^{-1}\L_{\t}=supp(\s_{\t})$. The moduli space of rank-two semistable vector-bundles on
$C$ with trivial determinant has dimension three, hence a theorem of Raynaud~[R] shows that
if $\t_1,\ldots,\t_4$ are generic line-bundles of degree one, there is no semistable
vector-bundle $V$ with trivial determinant such that
$$h^0(V\ot\t_i)>0,\qquad i=1,\ldots,4.$$
On the other hand, if $V$ is a degree zero rank-two sheaf on $C$ which is either singular or
locally-free non-semistable, then $h^0(V\ot\t)>0$ for any choice of a degree-one line-bundle
$\t$ on $C$. So let  $\t_1,\ldots,\t_4$ be generic, and set
$\L:=\L_{\t_1}\cap\cdots\cap\L_{\t_4}$; then 
$$\phi^{-1}\L=
\{[F]\in\cM_4|\ \hb{$F|_C$ is singular or locally-free non-semistable}\}.\leqno(1.3.2)$$ 
Let $\Si_C,B_C\ss\cM_4$ be defined as
$$\Si_C:=\{[F]\in\Si|\ sing [F]\cap C\not=\es\}\quad 
B_C:=\{[F]\in B|\ sing [F]\cap C\not=\es\},$$
where $\Si\ss\cM_4$ is the locus parametrizing strictly semistable (i.e.~non stable)
sheaves. (Thus $\Si_C\ss B_C$.) Let $V_C^0\ss\cM_4$ be given by
$$V_C^0:=\{[F]\in \cM_4^{lf}|\ \hb{$F|_C$ is not semistable}\},$$
and let $V_C$ be the closure of $V_C^0$. By Proposition~(1.1.1) we can rewrite~(1.3.2) as
$$\phi^{-1}\L=V_C\cup B_C.$$
Let $\tV_C,\tB_C\ss\cMt_4$ be the proper transforms of $V_C$ and $B_C$
respectively, and set $\Sit_C:=\pit^{-1}\Si_C$. The equality above gives
$$\phit^{-1}\L=\tV_C\cup\Sit_C\cup\tB_C,$$
hence by Corollary~(1.2.5) we get the following result.

\proclaim (1.3.3) Proposition.
The map 
$$H^q(\cMt_4;\ZZ)\to H^q(\tV_C\cup\Sit_C\cup\tB_C;\ZZ)$$
induced by inclusion is an isomorphism for $q\le 5$. 

\bsk
\n
{\bf 1.4. Elementary modifications.}
\bsk
\n
Let $[F]\in V_C^0$, and let
$$0\to L\to F|_C\to L^{-1}\to 0\leqno(1.4.1)$$
be the desemistabilizing sequence, i.e.~$L$ is  a line bundle of degree $d>0$. Let $\i\cl C\hra
X$ be inclusion, and let $G$ be the elementary modification of $F$ determined by~(1.4.1),
i.e.~the sheaf fitting into the exact sequence
$$0\to G\brel\a\over\to F\to \i_{*}L^{-1}\to 0.\leqno(1.4.2)$$
Then $G$ is locally-free of rank two, and $c_1(G)=-H$, $c_2(G)=(4-d)$. (First compute the
Chern classes of $\i_{*}L^{-1}$ by applying Grothendieck-Riemann-Roch to $\i$.) 

\proclaim (1.4.3) Lemma.
Keep notation as above. The vector-bundle $G$ is slope-stable.

\pf
Assume $G$ is not slope-stable. Since $\Pic(X)=\ZZ[H]$, there exists an injection
$\cO_X(k)\hra G$ with $k\ge 0$. Composing with $\a$ we get an injection $\cO_X(k)\hra F$,
contradicting (Gieseker-Maruyama) semistability of $F$.
\qed

\proclaim (1.4.4) Corollary.
Keeping notation as above, either $d=1$ or $d=2$.

\pf
By Serre duality    $h^2(G^{*}\ot G)=h^0(G^{*}\ot G)$. Since $G$ is
slope-stable $h^0(G^{*}\ot G)=1$,  hence $\chi(G^{*}\ot G)\le 2$.  Hirzebruch-Riemann-Roch
gives $\chi(G^{*}\ot G)=(4d-6)$, thus $d\le 2$.
\qed
\msk
By the above corollary we have a decomposition into locally
closed subsets
$$V_C^0=V_C^0(1)\coprod V_C^0(2),$$
where 
$$V_C^0(d):=\{[F]\in V_C^0|\ \hb{$F|_C\to L^{-1}\to 0$, $L$ a degree-$d$ line-bundle}\}.$$
We will describe $V_C^0(d)$ in terms of $\cM(H,4-d)$,  where 
$\cM(H,4-d)$ is the moduli space of semistable rank-two torsion-free sheaves on $X$ with
$c_1=H$, $c_2=(4-d)$. Let $\cM(H,4-d)^{lf}\ss\cM(H,4-d)$ be the subset parametrizing
locally-free sheaves. Keeping notation as above, let $E:=G(1)$; then $c_1(E)=H$,
$c_2(E)=(4-d)$, so that by~(1.4.3) $[E]\in\cM(H,4-d)^{lf}$. In order to reconstruct $F$ from
$E$, we notice that the bundle $E|_C$ comes with a canonical rank-one subsheaf: in fact
the long exact sequence of $Tor(\cdot,\cO_C)$ 
associated to~(1.4.2) gives  
$$0\to L^{-1}\ot\cO_C(-C)\brel\b\over\to G|_C\to L\to 0,$$
and if we tensor the above sequence with $\cO_C(1)$ we obtain
$$0\to L^{-1}\to E|_C\to L\ot K_C\to 0.\leqno(1.4.5)$$
One recovers $F$ from $E$ and~(1.4.5) as follows. First notice that $F(-C)\hra G$,
secondly that the restriction to $C$ of this inclusion has image equal to $\Im\b$, thus we
have an exact sequence
$$0\to F(-1)\to G\to \i_{*}L\to 0.$$  
Tensoring with $\cO_X(1)$ we see that $F$ is the elementary modification of $E$
associated to~(1.4.5). The following result says that  we can ``invert'' this construction.  

\proclaim (1.4.6) Lemma.
Let $d>0$, and let $[E]\in\cM^{lf}(H,4-d)$. Assume $E|_C$ fits into Exact
Sequence~(1.4.5), where $L$ is a line-bundle with $\deg L=d$. Let $F$ be the elementary
modification of $E$ associated to~(1.4.5), i.e.~we have
$$0\to F\brel\g\over\to E\to\i_{*}(L\ot K_C)\to 0.\leqno(1.4.7)$$
Then $F$ is a rank-two slope-stable vector-bundle with $c_1(F)=0$, $c_2(F)=4$, and $[F]\in
V_C^0(d)$. 

\pf
The Chern classes of $F$ are easily computed from the exact sequence defining $F$.
Furthermore, applying the functor $Tor(\cdot,\cO_C)$ to~(1.4.7) we get Exact
sequence~(1.4.1), so all we have to prove is that $F$ is slope-stable. Suppose $F$ is not
slope-stable.  Since $\Pic(X)=\ZZ[H]$, there is an injection $\cO_X(k)\hra F$, with $k\ge 0$.
Composing with $\g$ we get $\cO_X(k)\hra E$; since $E$ is semistable we must have $k=0$,
i.e.~we have a non-zero section $\s\in H^0(E)$. The restriction of $\s$ to $C$ is a section of
$L^{-1}$, hence zero because $\deg(L^{-1})<0$. Since $C\in|\cO_X(1)|$ the section $\s$ gives
rise to an injection $\cO_X(1)\hra E$, contradicting semistability of $E$. 
\qed
\msk
Let $\cE_d$ be a tautological vector-bundle on $X\tm \cM^{lf}(H,4-d)$ (it exists
by~[M2,(A.7)]). Let $\cQ_C(d)$ be the relative quot-scheme of $\cE_d|_{C\tm\cM^{lf}(H,4-d)}$
over $\cM^{lf}(H,4-d)$ parametrizing quotients $E|_C\to \xi$, where $\xi$ is a rank-one sheaf
of degree $(2+d)$, and let $\cQ_C^0(d)\ss\cQ_C(d)$ be the open
subset parametrizing locally-free quotients. Lemmas~(1.4.3)-(1.4.6) prove the following
result.

\proclaim (1.4.8).
{\rm There is an isomorphism between $\cQ_C^0(d)$ and $V_C^0(d)$, defined as follows: to a
quotient $E|_C\to L\ot K_C$, where $[E]\in\cM(H,4-d)^{lf}$ and $L$ is a degree-$d$ line-bundle,
we associate the point $[F]\in V_C^0(d)$, where $F$ is the locally-free sheaf fitting into
the exact sequence
$$0\to F\to E\to\i_{*}L\ot K_C\to 0.$$}
   
\bsk
\n
{\bf 1.5. $V_C^0(2)$  is nowhere dense.}
\bsk
\n
We will prove the following.

\proclaim (1.5.1) Proposition.
$V_C^0(1)$ is (open) dense in $V_C^0$.

The proof of the proposition will be given at the end of this subsection. First we describe
$\cM(H,2)$. Let $U:=T_{\PP^2}(-1)$;
this a slope-stable rank-two vector-bundle on $\PP^2$ with $c_1(U)=c_1(\cO_{\PP^2}(1))$,
$c_2(U)=1$.  (That it is slope-stable follows from the fact that $U$ has sections with an
isolated zero.) Let 
$$\psi\cl X\to |\cO_X(1)|^{*}\cong\PP^2$$ 
be the two-to-one branched cover, and set $W:=\psi^{*}U$. 

\proclaim (1.5.2) Lemma.
Keeping notation as above, $W$ is a slope-stable rank-two vector-bundle on $X$ 
with $c_1(W)=H$, $c_2(W)=2$. We have
$$W|_C\cong\cO_C\op K_C.\leqno(1.5.3)$$
Furthermore $\cM(H,2)=[W]$.  

\pf
Let's prove $W$ is stable. Since $U$ has
sections with isolated zeroes, so does $W$, hence we have an exact sequence
$$0\to \cO_X\to W\to I_Z(1)\to 0,$$
where $Z\ss X$ is a zero-dimensional subscheme of length $2$.
Since $\Pic(X)=\ZZ[H]$, it follows that $W$ is slope-stable. 
Now let $R:=\psi(C)$, a line in $\PP^2$. To prove~(1.5.3) it suffices to verify
that 
$$U|_R\cong\cO_R\op \cO_R(1).$$
This follows immediately from the exact sequence
$$0\to T_R\to T_{\PP^2}|_R\to N_{R/\PP^2}\to 0.$$
Let's prove the last statement. Since the expected dimension of $\cM(H,2)$   is zero, a
result of Mukai~[M2,(3.6)] gives that $\cM(H,2)$ consists of a single point,
hence $\cM(H,2)=[W]$. 
\qed

\proclaim (1.5.4) Corollary.
Keeping notation as above, $\dim V_C^0(2)=5$.

\pf
According to~(1.4.8) we have an identification $V_C^0(2)\cong\cQ^0_C(2)$. 
Thus we have a morphism 
$$\matrix{ V_C^0(2) & \to & \Pic^2(C) \cr
\left(W|_C\to L\ot K_C\to 0\right) & \mapsto & [L] \cr}.$$
The fiber over $[L]$ is the open subset of $\PP H^0(L\ot W|_C)$ corresponding to sections
with no zeroes. It follows from~(1.5.3) that
$$\dim\PP H^0(L\ot W|_C)=\cases{ 3 & if $L\not\cong K_C$, \cr
4 & if $L\cong K_C$. \cr}$$
The corollary follows immediately. 
\qed
\msk
\n
{\bf Proof of Proposition~(1.5.1).}
\hskip 2mm
By an argument similar to that proving Proposition~(1.13) of~[O3], one shows that every
irreducible component of $V_C^0$ has codimension at most $g(C)+1=3$. Since $\cM_4$ is of
pure dimension $10$, every component of $V_C^0$ has dimension at least $7$. The proposition
follows from Corollary (1.5.4).    
\bsk
\n
{\bf 1.6. Analysis of $\cM(H,3)$.}
\bsk
\n
We are mainly interested in the restriction to $C$ of vector-bundles parametrized by
$\cM(H,3)^{lf}$: the key results are stated in Propositions~(1.6.4)-(1.6.6). By~[M1] the moduli
space $\cM(H,3)$ is smooth symplectic of dimension $4$, and by~[O1] it is irreducible. 
Our first goal is to define a morphism 
$$\rho\cl\cM(H,3)\to |\cO_X(1)|.\leqno(1.6.1)$$

\proclaim (1.6.2) Lemma.
If $[E]\in\cM(H,3)$ then $h^0(E)=2$.

\pf
By Riemann-Roch $\chi(E)\ge 2$. By Serre duality $h^2(E)=\dim \Hom(E, \cO_X)$,
hence stability  gives $h^2(E)=0$. Thus $h^0(E)\ge 2$. Let $\s\in H^0(E)$ be non-zero. We
claim the quotient $Q:=E/\cO_X\s$ is torsion-free. Suppose $Q$ has torsion: then $\s$ must
vanish on a divisor, and since $\Pic(X)=\ZZ[H]$ this implies we have an injection
$I_W(k)\hra E$, where $k\ge 1$, and $W\ss X$ is a zero-dimensional subscheme. This  
contradicts semistability of $E$. Thus we have an exact sequence
$$0\to \cO_X\to E\to I_Z(1)\to 0,$$
where $Z\ss X$ is a zero-dimensional subscheme of length $c_2(E)=3$. Since $H\cdot H=2$,
we have $h^0(I_Z(1))\le 1$, hence $h^0(E)\le 2$.
\qed
\msk
In order to define $\rho$ choose a basis $\{\s,\tau\}$ of $h^0(E)$. We claim
$\s\wedge\tau\not\equiv 0$. Assume the contrary: then $\s$, $\tau$ generate a rank-one
subsheaf of $E$ with two linearly independent global sections, i.e.~a sheaf isomorphic to
$I_W(k)$, where $k\ge 1$ and $W\ss X$ is a zero-dimensional subscheme. This contradicts
semistability of $E$. Thus 
$$0\not=\s\wedge\tau\in H^0(\det E)=H^0(\cO_X(1)).$$
We set
$$\rho([E]):=(\s\wedge\tau)\in |\cO_X(1)|.$$

\proclaim (1.6.3) Lemma.
For all $D\in|\cO_X(1)|$, the fiber  $\rho^{-1}D$ has pure dimension $2$. In particular 
the map $\rho$ is surjective. 

\pf
Since $\dim\cM(H,3)=4$, every component  of $\rho^{-1}D$ has dimension at least $2$. It
follows from~[O2,(2-9)] that the symplectic form on $\cM(H,3)$ is identically zero on
$\rho^{-1}D$, hence every component of $\rho^{-1}D$ has dimension at most $2$. 
\qed

\proclaim (1.6.4) Proposition.
Let $[E]\in\cM(H,3)^{lf}$. Then $E|_C$ is not semistable if and only if $\rho([E])=C$. In this
case
$$E|_C\cong\xi\op (K_C\ot\xi^{-1}),$$
where $\xi$ is degree $3$ line-bundle.

\pf
Assume $\rho([E])=C$. Choose a basis $\{\s,\tau\}$ of $H^0(E)$, and  consider the
exact sequence 
$$0\to\cO_X^{(2)}\brel (\s,\tau)\over\lra E\to\i_*\eta\to 0.$$
Restricting to $C$ we get
$$0\to\xi\to E|_C\to \eta\to 0,\leqno(1.6.5)$$
where $\xi$, $\eta$ are rank-one sheaves. Thus $\xi$ is locally-free of rank-one, and a Chern
class computation gives $\deg\xi=3$. Since $C=(\s\wedge\tau)$ is smooth, $\s$ and
$\tau$ have no zeroes in common, and this implies $\eta$ is also locally-free. Since $\det
(E|_C)\cong K_C$, we have $\eta\cong\xi^{-1}\ot K_C$. Exact Sequence~(1.6.5) splits because
$H^1(\xi^{\ot 2}\ot K_C^{-1})=0$. Now suppose $E|_C$ is not semistable, and let~(1.6.5) be the
desemistabilizing sequence, i.e.~$\eta$ is a line-bundle and $d=\deg\xi\ge 2$. Let $F$ be the
elementary modification of $E$ defined by~(1.6.5), i.e.~we have 
$$0\to F\to E\to\i_*\eta\to 0.$$
Then $\rk F=2$, $c_1(F)=0$ and $c_2(F)=(3-d)\le 1$. Arguing as in the proof of Lemma~(1.4.3)
we see that $F$ is slope-semistable. Since there are no slope-semistable rank-two
vector-bundles on $X$ with $c_1=0$, $c_2=1$ (see the proof of~(1.1.1)), we conclude that
$c_2(F)=0$, hence $F\cong\cO_X^{(2)}$. This means that $\rho([E])=C$.
\qed

\proclaim (1.6.6) Proposition.
Let $[E]\in\cM(H,3)^{lf}$. Assume $D=\rho([E])$ is not equal to $C$, and set $D\cdot
C=p+p'$. There is an exact sequence 
$$0\to\cO_C(p)\to E|_C\to\cO_C(p')\to 0,\leqno(1.6.7)$$
which is split if $p\not=p'$. 

\pf
Choose a basis $\{\s,\tau\}$ of $H^0(E)$. Since $\s(p)$, $\tau(p)$
are linearly dependent, there exists a non-trivial linear combination $\l\s+\mu\tau$
which is zero at $p$. Let
$$\e:=(\l\s+\mu\tau)|_C\in H^0(E|_C).$$
Since $\e$ vanishes at $p$, it defines a non-zero map $\cO_C(p)\to E|_C$. 
By Proposition~(1.6.4) the bundle $E|_C$ is semistable, hence $\e$ vanishes only at
$p$ and with multiplicity one. This gives Exact Sequence~(1.6.7).  Reversing the roles of
$p$ and $p'$, we see that if $p\not=p'$ the sequence is split.
\qed    

\proclaim (1.6.8) Remark.
{\rm The proof of Proposition~(1.6.4) gives an isomorphism
$$\matrix{\rho^{-1}C & \brel\sim\over\to & \Pic^3(C) \cr
[E] & \mapsto & [\xi] \cr}.$$
Obviously this holds also if we replace $C$ by any smooth $D\in |\cO_X(1)|$. 
For a  general $D$, the fiber $\rho^{-1}D$ is identified with the moduli space of rank-one
torsion-free sheaves of degree $3$ on $D$. The identification is given by associating to $[E]$
the subsheaf of $E|_D$ generated by $H^0(E)$.}

\bsk
\n
{\bf 1.7. Birational description of $V_C$.}
\bsk
\n
Let $t\in\cQ_C(1)$, and let
$$E|_C\brel f\over\to\xi\to 0$$
be the quotient represented by $t$: thus $[E]\in\cM(H,3)^{lf}$, and $\xi$ is a 
rank-one sheaf on $C$ of degree $3$. Since $E|_C$ is locally-free of
degree $2$, the kernel of $f$ is a rank-one locally-free sheaf of degree
$(-1)$, say $L^{-1}$. Thus we can define a morphism
$$\matrix{ \cQ_C(1) & \brel\vf\over\lra & \Pic^1(C)\tm\cM(H,3)^{lf}
\cr t & \mapsto & ([L],[E]) \cr}.$$
We let $\vf_0$ be the restriction of $\vf$ to $\cQ_C^0(1)$. Clearly we
have isomorphisms
$$\leqalignno{ \vf^{-1}([L],[E]) \cong & \PP H^0(L\ot E|_C), & (1.7.1) \cr
\vf_0^{-1}([L],[E]) \cong & 
\PP\{\s\in H^0(L\ot E|_C)|\ \ \hb{$\s$ has no zeroes}\}. & (1.7.2) \cr}$$
Let's examine $\vf$ over the open subset
$\cU\ss\Pic^1(C)\tm\cM(H,3)^{lf}$ consisting of couples $([L],[E])$ such
that $h^0(L)=0$ and $\rho([E])$ intersects $C$ in exactly two 
points. Set
$$\eqalign{ P_C:= & \vf^{-1}(\cU), \cr
P_C^0:= & \vf_0^{-1}(\cU). \cr}$$

\proclaim (1.7.3) Claim.
Both $P_C$ and $P_C^0$ are fibrations over $\cU$, the former with fibers
$\PP^1$, the latter with fibers $\CC^{*}$.
In particular $\dim P_C=\dim P_C^0=7$.

\pf
Let $([L],[E])\in\cU$, and let $\rho([E])\cap C=\{p,p'\}$. By~(1.6.6)
$$\leqalignno{ E|_C\cong & \cO_C(p)\op\cO_C(p'), & (1.7.4)\cr
H^0(L\ot E|_C)= & H^0(L(p))\op H^0(L(p')). & (1.7.5)\cr}$$
Since $K_C\sim p+p'$, and since $L\not\cong[p']$, we see that
$h^0(L(p))=1$. Similarly $h^0(L(p'))=1$, hence the space of sections
in~(1.7.5) is two-dimensional. By~(1.7.1) we conclude that the fibers
of $\vf$ are isomorphic to $\PP^1$. An easy argument identifies $P_C$ 
with the projectivization of a direct image sheaf over $\cU$ with fiber
$H^0(L\ot E|_C)$ over the point $([L],[E])$. This proves the result for
$\vf$. To finish the proof we remark that, since $h^0(L)=0$, a section
$$\s=(\s_1,\s_2)\in H^0(L(p))\op H^0(L(p'))$$
has a zero if and only if $\s_1=0$ or $\s_2=0$. 
\qed
\msk
The following result will allow us to forget about the complement of
$P_C^0$ in $\cQ_C^0(1)$.

\proclaim (1.7.6) Proposition.
Keeping notation as above, $P_C^0$ is dense in $\cQ_C^0(1)$.

\pf
This is a dimension count. Stratify the complement of $\cU$ in
$\Pic^1(C)\tm\cM(H,3)^{lf}$ according to the dimension of the fibers of
$\vf_0$ (see~(1.7.2)). Using~(1.6.4)-(1.6.6), one easily verifies that for
each stratum $\cS$,
$$\dim\vf_0^{-1}(\cS)<7.$$
On the other hand (see the proof of~(1.5.1)) every irreducible component of $V_C^0(1)$ has
dimension at least $7$. By~(1.4.8) $\cQ_C^0(1)\cong V_C^0(1)$, hence the
above inequality shows that $\vf_0^{-1}(\cS)$ is not dense in any
component of $\cQ_C^0(1)$.
\qed
\msk
Restricting the isomorphism $\cQ_C^0(1)\brel\sim\over\to V_C^0(1)$ to
$P_C^0$, we get a map
$$\psi_0\cl P_C^0\to V_C,$$
 which is an isomorphism onto its image. The following result is an
immediate consequence of~(1.5.1) and~(1.7.6). 

\proclaim (1.7.7) Corollary.
The image $\psi_0(P_C^0)$ is dense in $V_C$. In particular $V_C$ is
irreducible of dimension $7$. 

\bsk
\n
{\bf 1.8. Extension of $\psi_0$.}
\bsk
\n
We will extend $\psi_0$ to a map 
$$\psi\cl P_C\to V_C\leqno(1.8.1)$$
Let $t\in P_C$
correspond to the quotient
$$E|_C\to\xi\to 0,$$
and let $F$ be the associated elementary modification, i.e.~we have
$$0\to F\to E\to\i_{*}\xi\to 0.$$
The sheaf $F$ is torsion-free of rank two, with $c_1(F)=0$, $c_2(F)=4$.
Since the family of quotient sheaves $\xi$ form a flat family over $P_C$,
the sheaves $F$ fit together to give a shef $\cF$ on $X\tm P_C$, flat over
$P_C$. Let 
$$N_C:=P_C\sm P_C^0.$$

\proclaim (1.8.2) Claim.
Let $t\in N_C$, and let
$$F=\cF|_{X\tm\{t\}}.$$
Then $F$ is  a strictly (Gieseker-Maruyama) semistable sheaf. More
precisely there is an exact sequence
$$0\to I_Z\to F\to I_W\to 0,\leqno(1.8.3)$$
where $Z,W\ss X$ are zero-dimensional subschemes of length $2$, with
$Z\not= W$. 

\pf
Let $\vf(t)=([L],[E])$. By~(1.6.6)
$$E|_C\cong\cO_C(p)\op\cO_C(p'),$$
where $p\not=p'$. We can assume (see the proof of~(1.7.3)) that $t$
corresponds to the quotient of the inclusion
$$L^{-1}\brel\a\over\hra\cO_C(p)\hra E|_C.$$
The map $\a$ vanishes on a divisor $Z\ss C$ of degree $2$. There is a
global section $\tau\cl\cO_X\to E$ which restricted to $C$ gives the
non-zero section of $\cO_C(p)$ (see the proof of~(1.6.6)), hence
$$\tau(I_Z)|_C=\Im\a.$$
Thus $\tau(I_Z)$ is a subsheaf of the elementary modification $F$, i.e.~we
have an exact sequence
$$0\to I_Z\to F\to \eta\to 0,$$
for some rank-one sheaf $\eta$. We claim $\eta$ is torsion-free. First, a
local computation shows that $\eta$ is locally-free at points of $C$.
Secondly, since outside $C$ the sheaves $F$ and $E$ are isomorphic, also
$\eta$ and $E/\Im\tau$ are isomorphic outside $C$. By slope-stability,
the section $\tau$ has isolated zeroes, hence $E/\Im\tau$ is
torsion-free, and we conclude that $\eta$ is torsion-free outside $C$.
Thus $\eta\cong I_W(k)$ for some integer $k$. Since $c_1(F)=0$, we have
$k=0$, and since $c_2(F)=4$ the length of $W$ is $2$. This proves $F$
fits into Exact Sequence~(1.8.3). Finally, since $\eta$ is locally-free on
$C$, $W$ does not intersect $C$. Since $Z\ss C$, we conclude $Z\not=W$.
\qed
\msk
By the above claim the sheaves $\cF|_{X\tm\{t\}}$ are semistable for all
$t\in P_C$, hence $\cF$ induces a morphism $\psi\cl P_C\to\cM_4$,
which is equal to $\psi_0$ on $P_C^0$. Since $\Im\psi_0\ss V_C$, and
$V_C$ is closed, we have $\Im\psi\ss V_C$. Thus we have defined the
desired extension~(1.8.1) of $\psi_0$. By Claim~(1.8.2) we have
$$\psi\cl P_C\to V_C\sm\O,$$
where we recall that $\O$ parametrizes strictly semistable sheaves of
the form $I_Z\op I_Z$. Let's show that $\psi$ lifts to a map
$$\wt{\psi}\cl P_C\to \tV_C.$$
Outside $\O$ the desingularization map $\pit$ is the blow-up of
$(\Si\sm\O)$. Since $\psi^{-1}(\Si)=N_C$, and since $N_C$ is a
divisor in the smooth variety $P_C$, hence Cartier, $\psi$ lifts to
$\cMt_4$ by the universal property of blow-up. Obviously $\wt{\psi}(P_C)$
is contained in $\tV_C$, and is dense in it. We choose a smooth projective
completion $V'_C$ of $P_C$ such that $\wt{\psi}$ extends to a morphism
$$\Psi\cl V'_C\to \tV_C.$$
Summing up, we have proved the following result.

\proclaim (1.8.4) Proposition.
Keep notation as above. Then
\msk
\item{\rm 1.}
$\Psi\cl V'_C\to \tV_C$ is a birational morphism,
\item{\rm 2.}
$\Psi(P_C)$ is dense in $\tV_C$, in particular $\tV_C$ is irreducible,
\item{\rm 3.}
$\Psi(N_C)\ss\Sit_C$, in particular $\tV_C\cap\Sit_C\not=\es$.

A key observation for the proof of
$h^{2,0}(\cMt_4)=1$ is the following.

\proclaim (1.8.5) Lemma.   
The restriction map
$$H^0(\O^2_{V'_C})\to H^0(\O^2_{N_C})$$
is injective.

\pf
Let $\tau$ be a regular two-form on $V'_C$ such that $\tau|_{N_C}\equiv
0$;
 $V'_C$ is birational to $\PP^1\tm \Pic^1(C)\tm\cM(H,3)$, hence there
exists a regular two-form $\e$ on  $\Pic^1(C)\tm\cM(H,3)$ such that  
$$\tau|_{P_C}=\vf^{*}\e.$$
Since the restriction of $\vf$ to $N_C$ is a two-to-one covering of  
$\Pic^1(C)\tm\cM(H,3)^{lf}$, we get that $\e$ is identically zero on 
$\Pic^1(C)\tm\cM(H,3)^{lf}$. But $\cM(H,3)^{lf}$ is dense in $\cM(H,3)$,
hence $\e\equiv 0$. Thus $\tau\equiv 0$.
\qed 
\bsk
\n
{\bf 1.9. Analysis of $\tB_C$ and $\Sit_C$.}
\bsk
\n
Consider the restriction of the determinant map $\phi$ (see~(1.1)) to $B$:
$$\matrix{B & \brel{\phi|_{B}}\over\lra & X^{(4)}\cr
[F] & \mapsto & sing [F].\cr}$$
Let $X^{(4)}_{(4)}\ss X^{(4)}$ be the open stratum of Stratification~(1.2.2),
i.e. the partition is $(4):=(4,0,0,0)$. Set
$$\eqalign{ B^0:= & \phi^{-1}X^{(4)}_{(4)},\cr
\Si^0:= & \Si\cap B^0.\cr}$$

\proclaim (1.9.1) Lemma.
The map $\phi|_{B^0}\cl B^0\to X^{(4)}_{(4)}$ is a
$\PP^1$-fibration in the analytic topology. In particular $B^0$ is smooth of
dimension $9$. 
Furthermore $\Si^0$ is a smooth divisor in $B^0$. 

\pf
Let $\g\in X^{(4)}_{(4)}$, so $\g=q_1+\cdots+q_4$, where the $q_j$ are
pairwise distinct. Let $Quot(\cO_X^{(2)},\g)$ be the quot-scheme
parametrizing quotients  
$$0\to F\to\cO_X^{(2)}\to \bigoplus\limits_{j=1}\limits^4\CC_{q_j}\to
0.$$ 
Thus
$$Quot(\cO_X^{(2)},\g)\cong\PP^1\tm\PP^1\tm\PP^1\tm\PP^1.$$
Every sheaf parametrized by
$\phi^{-1}\g$ is isomorphic to the kernel $F$ of some quotient in 
$Quot(\cO_X^{(2)},\g)$. The group $\PP
{\rm Aut}(\cO_X^{(2)})\cong\PGL(2)$ acts on this quot-scheme: this is the
diagonal action on $\PP^1\tm\cdots\tm\PP^1$.   As is easily verified
semistability of a sheaf is equivalent to semistability of the
corresponding point in $\PP^1\tm\cdots\tm\PP^1$ with respect to the
symmetric linearization. Hence
$$\phi^{-1}\g=
\PP^1\tm\PP^1\tm\PP^1\tm\PP^1//\PGL(2)\cong\PP^1,$$
so that the fibers of $\phi|_{B^0}$ are $\PP^1$'s. It
is clear that $\phi|_{B^0}$ is locally trivial in the analytic topology.
Finally, let $\g\in X^{(4)}_{(4)}$; since $\phi^{-1}\g\cap\Si$
consists of the three points representing strictly semistable orbits,
$\Si\cap B^0$ is a smooth divisor in $B^0$. 
\qed
\msk
Let
$$\eqalign{\cV_C:= & \{\g\in X^{(4)}_{(4)}|\ \ \#\g\cap C=1\},\cr
B^0_C:= & \phi^{-1}\cV_C. \cr}$$
Clearly $B^0_C\ss B_C$.

\proclaim (1.9.2) Corollary.
\msk
\item{\rm 1.}
The map $\phi|_{B_C^0}\cl B_C^0\to \cV_C$ is a  $\PP^1$-fibration,
locally trivial in the analytic topolgy. 
\item{\rm 2.}
$B_C^0$ is a smooth irreducible locally closed codimension one subset of
$B^0$. In particular $\dim B_C^0=8$.  
\item{\rm 3.}
$\Si^0$ and $B_C^0$ intersect transversely (inside $B^0$),
and the intersection is given by
$$\Si_C^0:=
\{[I_Z\op I_W]|\ \ \hb{$Z$, $W$ reduced disjoint, $\# Z\cap C=1$,
$W\cap C=\es$}\}.$$

\pf
Item~(1) follows from Lemma~(1.9.1). Items~(2)-(3) are easily verified.
\qed
\msk
Let $\tB_C^0\ss\cMt_4$ be the strict transform of $B_C^0$; clearly
$\tB_C^0\ss\tB_C$. 

\proclaim (1.9.3) Proposition.
\msk
\item{\rm 1.}
$\tB_C^0$ is open dense in $\tB_C$.
\item{\rm 2.} 
The map $\tB_C^0\brel\phit\over\to \cV_C$ is a
$\PP^1$-fibration, locally trivial in the analytic topology.  
\item{\rm 3.}
$\tB_C^0$ is smooth irreducible of dimension $8$.

\pf
Let's prove Item~(1). First we claim that every irreducible component of
$\tB_C$ has dimension at least $8$. This is equivalent to showing that
every component of $B_C$ has dimension at least $8$. Let 
$$X_C^{(4)}:=
\{\g\in X^{(4)}|\ \ {\rm supp}\g\cap C\not=\es\},$$ 
so that  
$$B_C=\left(\phi|_{B}\right)^{-1}X_C^{(4)}\leqno(1.9.4).$$
Since  $X^{(4)}$ is the Quotient of a smooth variety by a finite group,
$X_C^{(4)}$ is a $\QQ$-Cartier divisor in $X^{(4)}$; by~(1.9.4) every
component of $B_C$ has codimension at most one in $B$, and this
proves our claim.  Now consider the stratification of $X_C^{(4)}$ induced by
Stratification~(1.2.2). If $\cS_p\ss X_C^{(4)}$ is the stratum indicized by
the partition $p=(p_1,\ldots,p_4)$, then by Corollary~(1.2.5) 
$$\dim\phit^{-1}\cS_p\le\dim\cS_p+5-\sum\limits_{j=1}\limits^4 p_j.$$
Since $\dim\cS_p=(2\sum\limits_{i=1}\limits^{4}p_i-1)$, we
get 
$$\dim\phit^{-1}\cS_p\le \sum\limits_{i=1}\limits^{4}p_i+4.$$
The quantity on the right is strictly less than $8$, unless $p=(4)$ (in which
case it equals $8$). Thus $\phit^{-1}\cS_{(4)}$ is dense in $\tB_C$. Since
$\cV_C$ is dense in $\cS_{(4)}$, this proves Item~(1). 
To prove Item~(2) notice that $B_C^0$ does not intersect $\O$, hence
$\tB_C^0$ is the blow-up of $B_C^0$ at the (scheme-theoretic)
intersection $B_C^0\cap\Si$.  By Item~(3) of Corollary~(1.9.2) the
intersection is a Cartier divisor, hence $\tB_C^0\cong B_C^0$. Thus
Item~(2) follows from Item~(1) of Corollary~(1.9.2).  Item~(3) follows at
once from Item~(2). 
\qed
\msk
Similar results hold for $\Sit_C$. 

\proclaim (1.9.5) Proposition.
$\Sit_C$ is irreducible of dimension $8$. 

\pf
Consider Kirwan's desingularization $\pih\cl\cMh_4\to\cM_4$
(see~[O4,(1.8)]), and let $\Sih_C:=\pih^{-1}\Si_C$. Let
$\a\cl\cMh_4\to\cMt_4$ be the contraction along $\Oh$ defined
in~[O4,(2)]. Since $\Sit_C=\a(\Sih_C)$, it suffices to prove $\Sih_C$ is
irreducible of dimension $8$. Letting $\Sih:=\pih^{-1}\Si$, it follows
from~[O4,(1.7)] that $\Sih$ is a $\PP^1$-fibration over 
$$\Si_R//\PGL(N)=Bl_{\O}(\Si).$$
Since the inverse image of $\Si_C$ under the blow-up
$Bl_{\O}(\Si)\to\Si$ is irreducible of dimension $7$, we conclude that 
$\Sih_C$ is irreducible of dimension $8$.
\qed
\msk
Let $\Sit_C^0:=\pit^{-1}\Si_C^0$. The following result is an immediate
consequence of Proposition~(1.9.5).

\proclaim (1.9.6) Corollary.
$\Sit_C^0$ is a smooth open dense subset of $\Sit_C$.

Let $S_C:=\Sit_C^0\cap\tB_C^0$.

\proclaim (1.9.7) Proposition.
$S_C$ is a section of the $\PP^1$-fibration $\Sit_C^0\to\Si_C^0$.

\pf
This is an immediate consequence of Item~(3) of Corollary~(1.9.2).
\qed
\msk
Let $\Si'_C$ be a desingularization of $\Sit_C$. Since by~(1.9.6),
$\Sit_C^0$ is a smooth open subset of $\Sit_C$, we can think
$\Sit_C^0\ss\Si'_C$, hence also $S_C\ss\Si'_C$. The following immediate
consequence of Proposition~(1.9.7) will be important in the proof that
$h^{2,0}(\cMt_4)=1$. 

\proclaim (1.9.8) Corollary.
Keeping notation as above, the restriction map
$$H^0(\O^2_{\Si'_C})\to H^0(\O^2_{S_C})$$
is injective.

\vfill
\eject
\n
{\bf 1.10. Proof that $b_0(\cMt_4)=h^{2,0}(\cMt_4)=1$.}
\bsk
\n
Let's prove $\cMt_4$ is connected. By Proposition~(1.3.3) it suffices to
show that
$$\tV_C\cup\Sit_C\cup\tB_C\leqno(1.10.1)$$
is connected. Each of $\tV_C,\Sit_C,\tB_C$ is irreducible, by~(1.8.4),
(1.9.5), and~(1.9.3) respectively. By Proposition~(1.8.4)
$\tV_C\cap\Sit_C\not=\es$, and clearly $\Sit_C\cap\tB_C\not=\es$.
Thus~(1.10.1) is connected. Now we prove that $h^{2,0}(\cMt_4)=1$.
Let 
$$\Psi\cl V'_C\to\tV_C$$
be as in~(1.8.4), and let 
$$f\cl \Si'_C\to\Sit_C\qquad g\cl B'_C\to\tB_C$$
be desingularizations of $\Sit_C$, $\tB_C$ respectively. By
Proposition~(1.3.3) the map
$$\matrix{ 
H^{2,0}(\cMt_4) & \to & 
H^{2,0}(V'_C) \op H^{2,0}(\Si'_C) \op H^{2,0}(B'_C) \cr 
[\tau] & \mapsto & (\Psi^{*}[\tau],
f^{*}[\tau],g^{*}[\tau])\cr}\leqno(1.10.2)$$ 
is injective. We claim that also $g^{*}$ is injective. So let $\tau$ be a
two-form on $\cMt_4$, and assume that $g^{*}\tau\equiv 0$. If $S_C$ is
as in~(1.9.7), then 
$$f^{*}\tau|_{S_C}\equiv g^{*}\tau|_{S_C}\equiv 0.$$
By Corollary~(1.9.8) we conclude that $f^{*}\tau\equiv 0$. Letting
$N_C\ss P_C$ be as in~(1.8.4), it follows from $f^{*}\tau\equiv 0$ that
$\Psi^{*}\tau|_{N_C}\equiv 0$. Thus Lemma~(1.8.5) gives that 
$\Psi^{*}\tau\equiv 0$. By injectivity of~(1.10.2) we get that
$\tau\equiv 0$. This proves $g^{*}$ is injective. We claim that
$h^{2,0}(B'_C)=1$. In fact, by Proposition~(1.9.3) there is a rational
dominant map 
$$B'_C\cdots >X^{[3]}$$
with generic fiber isomorphic to $\PP^1$, and since $h^{2,0}(X^{[3]})=1$ it
follows that  $h^{2,0}(B'_C)=1$. Since $g^{*}$ is injective we get that
$h^{2,0}(\cMt_4)\le 1$. But we know $h^{2,0}(\cMt_4)\ge 1$, hence 
$h^{2,0}(\cMt_4)=1$.
\bsk
\bsk
\n
{\bf 2.  $\cMt_4$ is simply-connected.} 
\bsk
\n
We will show that $\cMt_4$ is birational to an open subset $\cJ$ of the Jacobian fibration
parametrizing (stable) degree-six line-bundles on curves in $|\cO_X(2)|$.   This implies that
$\pi_1(\cJ)$ surjects onto $\pi_1(\cMt_4)$, hence it will suffice to show that $\cJ$ is
simply connected. The latter result is proved by a monodromy
argument.  
\bsk
\n
{\bf 2.1. The Jacobian fibration.}
\bsk
\n
Let $C\ss|\cO_X(2)|$ be a reduced irreducible curve, $\i\cl C\hra X$ be inclusion, and $L$ a
degree-six line-bundle on $C$. We let $\cN$ be Simpson's moduli space~[LP2,S] of pure
one-dimensional sheaves $\xi$ on $X$, stable with respect to $\cO_X(1)$,  with
$$\chi(\xi(n))=\chi(\i_{*}L(n)).$$
Thus $[\i_{*}L]$ is a typical point of $\cN$. There is a morphism
$$\matrix{\cN & \brel\rho\over\lra & |\cO_X(2)|\cong \PP^5 \cr
                   \xi & \mapsto & {\rm supp}\xi. \cr}\leqno(2.1.1)$$
By results of Mukai~[M1] we know that $\cN$ is smooth of dimension $10$.  We will be
interested in a certain open subset of $\cN$ defined as follows. First, let $U\ss |\cO_X(2)|$
be the open subset parametrizing reduced curves. The following result is a straightforward
application of the definition of stability according to Simpson~[LP2,S].

\proclaim (2.1.2) Lemma.
Let $[C]\in U$, and let $L$ be a line-bundle on $C$ such that:
\item{\rm 1.}
if $C$ is irreducible, the degree of $L$ is $6$,
\item{\rm 2.}
if $C=C_1\cup C_2$ and $L_i:=L|_{C_i}$, the degree of $L_i$ is $3$. 
\msk
\n
Then $\i_{*}L$ is a stable pure one-dimensional sheaf on $X$. 

\proclaim (2.1.3) Definition.
{\rm Let $\cJ\ss\cN$ be the open set parametrizing sheaves
$\i_{*}L$, where $[C]\in U$, $\i\cl C\hra X$ is inclusion, and $L$ is a line-bundle on $C$
satisfying the hypotheses of Lemma~(2.1.2).  We will often denote by $(C,L)$ the point
$[\i_{*}L]\in \cJ$. }

Since $\cJ$ is open in $\cN$, the results of Mukai mentioned above give the following.

\proclaim (2.1.4) Proposition.
$\cJ$ is smooth of dimension $10$.

\bsk
\n
{\bf 2.2. $\cMt_4$ is birational to $\cJ$.}
\bsk
\n
Notice that if $(C,L)\in\cJ$, then $\chi(L)=2$, hence $h^0(L)\ge 2$. Let
$$\cJ^0:=\{(C,L)\in\cJ|\ \ \hb{$C$ is smooth, $h^0(L)=2$ and $L$ is globally generated}\}.$$
As is easily verified $\cJ^0$ is open and dense in $\cJ$. Let $(C,L)\in\cJ^0$. Following
Lazarsfeld~[La] we will associate to $(C,L)$ a stable rank-two vector-bundle on $X$ with
$c_1=0$, $c_2=4$; this construction 
 will define an isomorphism between $\cJ^0$ and an open
subset of $\cMt_4$. Since $L$ is globally generated, the evaluation map
$H^0(L)\ot\cO_X\to\i_{*}L$ is surjective: let $E$ be the sheaf on $X$ fitting into the exact
sequence  
$$0\to E\brel\e\over\to H^0(L)\ot\cO_X\to\i_{*}L\to 0.\leqno(2.2.1)$$

\proclaim (2.2.2) Lemma.
Keeping notation as above, $E$ is a slope-stable rank-two vector bundle on $X$ with Chern
classes $c_1(E)=-2H$, $c_2(F)=6$.

\pf
The Chern classes are easily computed from~(2.2.1). To show stability, consider the exact
sequence 
$$0\to L^{-1}\to H^0(L)\ot\cO_C\to  L\to 0.\leqno(2.2.3)$$
By~(2.2.1) we have
$$\Im(\e|_C)=L^{-1}.\leqno(2.2.4)$$
Now suppose $E$ is
not stable: since $\Pic(X)=\ZZ [H]$ there is an injection of sheaves
$\cO_X(k)\brel\a\over\hra E$, where $k\ge -1$. By~(2.2.4)
$$(\e\circ\a)|_C\in H^0(L^{-1}(-k)).$$
Since $\deg L=6$, we have $\deg L^{-1}(-k)\le -2$. Thus $\e\circ\a$ vanishes on
$C$, i.e.
$$\e\circ\a\in H^0(\cO_X^{(2)}(-k)(-C))=H^0(\cO_X^{(2)}(-k-2)).$$
Since $k\ge -1$, this last group is zero, contradiction.
\qed
\msk
Set $F:=E(1)$; by the above lemma $F$ is a slope-stable rank-two vector-bundle on
$X$ with $c_1(F)=0$, $c_2(F)=4$.   Thus we can define $\Phi^0\cl\cJ^0\to\cM_4$ by setting
$$\matrix{\cJ^0 & \brel\Phi^0\over\lra  & \cM_4 \cr
                 (C,L) & \mapsto & [F]. \cr}$$
We will identify an open subset $\cM_4^0$ of $\cM_4$ such that $\Phi^0$ is an isomorphism
onto $\cM_4^0$. Let $[F]\in\cM_4^{lf}$, and set $G=F(1)$. Then $c_1(G)=2H$, $c_2(G)=6$; by
Riemann-Roch we get $\chi(G)=2$. By Serre duality $h^2(G)=h^0(G^{*})$, hence stability
gives $h^2(G)=0$. Thus $h^0(G)\ge 2$, therefore the locus where $h^0(G)=2$ is open in
$\cM_4$. If $h^0(G)=2$, consider the evaluation map
$$H^0(G)\ot\cO_X\brel f\over\to G.$$
Since $G$ is slope-stable and $\Pic(X)=\ZZ[H]$, the determinant of $f$ is not
identically zero, so $(\det f)\in |\cO_X(2)|$. We set
$$\cM_4^0:=\{[F]\in\cM_4^{lf}|\  \ \hb{$h^0(F(1))=2$ and $(\det f)$ is smooth}\}.$$
Clearly $\cM_4^0$ is open in $\cM_4$. 

\proclaim (2.2.5) Proposition.
The map $\Phi^0$ is an isomorphism onto $\cM_4^0$.

\pf
Let's show $\Phi^0(\cJ^0)\ss\cM_4^0$. Tensoring~(2.2.1) by $\cO_X(2)$ and observing that
$G=E(2)$, we get
$$0\to G\to \cO_X(2)^{(2)}\to \i_{*}L(2)\to 0.\leqno(2.2.6)$$
The long exact sequence of $Tor(\cdot,\cO_C)$ associated to~(2.2.6)
gives
$$0\to L\to G|_C\to L^{-1}(2)\to 0.$$
From~(2.2.6) we get a natural injection $\b\cl\cO_X(2)^{(2)}(-C)\hra G$, and
$\Im(\b|_C)=L$. Thus we have an exact sequence
$$0\to\cO_X^{(2)}\brel\g\over\to G\to\i_{*}L^{-1}(2)\to 0.\leqno(2.2.7)$$
By adjunction $\cO_C(2)\cong K_C$, hence Serre duality gives $h^0(L^{-1}(2))=h^1(L)$. By
hypothesis $h^0(L)=\chi(L)$, hence $h^1(L)=0$. Thus the cohomology long exact sequence
of~(2.2.7) gives $h^0(G)=2$. Furthermore the map $\g$ drops rank along $C$, which is
smooth, hence $[F]\in\cM_4^0$. Now we define an inverse
$$(\Phi^0)^{-1}\cl\cM_4^0\to\cJ^0.$$
Let $[F]\in\cM_4^0$, and set $G=F(1)$. Since $[F]\in\cM_4^0$, we have an exact sequence
$$0\to H^0(G)\ot\cO_X\brel f\over\to G\to \i_{*}\xi\to 0,\leqno(2.2.8)$$
where $(\det f)\in|\cO_X(2)|$ is smooth: set $C:=(\det f)$. Since $C$ is smooth, the map $f$
is nowhere zero, hence $\xi$ is a line-bundle on $C$: we can write $\xi=L^{-1}\ot K_C$,
where $L$ is a line-bundle. By a Chern class computation one gets $\deg L=6$. Since
$h^0(G)=2$, the cohomology  long exact sequence of~(2.2.8) gives $h^0(L^{-1}\ot K_C)=0$.
Thus $h^1(L)=0$, so that $h^0(L)=2$. Furthermore, the long exact sequence of
$Tor(\cdot,\cO_C)$ associated  to~(2.2.8) gives
$$0\to L\to G|_C\to L^{-1}\ot K_C\to 0.$$
Since $(\Im f)|_C=L$, the line-bundle $L$ is generated by global sections. We
have proved $(C,L)\in\cJ^0$. As is easily verified  the map
$$\matrix{\cM_4^0 & \to & \cJ^0 \cr
                  [F] & \to & (C,L) \cr}$$
is the inverse of $\Phi^0$.
\qed
\msk
We can view $\Phi^0$ as a map to $\cMt_4$, because $\cM_4^0$ is in the stable locus of
$\cM_4$. Since $\cMt_4$ is irreducible, Proposition~(2.2.5) implies that $\Phi^0$ extends
to a birational map
$$\Phi\cl \cJ\cdots >\cMt_4.$$
Let $I(\Phi)\ss\cJ$ be the indeterminacy locus of $\Phi$. By Proposition~(2.1.4) $\cJ$ is
smooth, hence $I(\Phi)$ has codimension at least two and the map induced by inclusion
$$\pi_1(\cJ\sm I(\Phi))\to \pi_1(\cJ)$$
is an isomorphism. On the other hand, since $\Phi$ is a birational map and $\cMt_4$ is
smooth, the map
$$\pi_1(\cJ\sm I(\Phi))\to\pi_1(\cMt_4)$$
is surjective. Thus we have proved the following result. 

\proclaim (2.2.9) Proposition.
The map $\Phi$ induces a surjection $\pi_1(\cJ)\to\pi_1(\cMt_4)$.

\bsk
\n
{\bf 2.3. $\cJ$ is simply-connected.}
\bsk
\n
Let $f\cl X\to|\cO_X(1)|^{*}\cong\PP^2$ be the two-to-one cover, and let $B\ss\PP^2$ be
the (sextic) branch curve. 
Let $V\ss U$ be the open subset parametrizing smooth curves $C\in|\cO_X(2)|$. Then
$$V=U\sm(\D\cap\L),\leqno(2.3.1)$$
where
$$\eqalign{\D:= & \{C\in|\cO_X(2)|\ |\ \ \hb{$f(C)$ is a smooth conic tangent to $B$}\}, \cr
 \L:= & \{C\in|\cO_X(2)|\ |\ \ \hb{$f(C)$ is a singular reduced conic}\}. \cr}$$
Let $\rho\cl\cJ\to U$ be the canonical map~(2.1.1), and set $\cJ_V:=\rho^{-1}V$. By
Proposition~(2.1.4) $\cJ$ is smooth, hence the map 
$$j_{\#}\cl\pi_1(\cJ_V)\to\pi_1(\cJ)$$
induced by inclusion $j\cl\cJ_V\hra\cJ$ is surjective. We will show $j_{\#}$ is
trivial: this will prove $\cJ$ is simply-connecetd. The map $\cJ_V\to V$ is a fibration with
fibers $5$-dimensional Jacobians. The homotopy exact sequence of a fibration gives an
exact sequence 
$$\pi_1(Jac(C))\to\pi_1(\cJ_V)\to\pi_1(V)\to \{1\},$$
 where $C\in|\cO_X(2)|$ is a fixed smooth curve. 

\proclaim (2.3.2) Lemma.
 Keeping notation as above, the restriction of $j_{\#}$ to $\pi_1(Jac(C))$ is trivial.

\pf
We have an isomorphism $\pi_1(Jac(C))\cong H_1(C;\ZZ)$. As is easily seen the vanishing
cycles on $C$ for the family of curves parametrized by $V$ generate all of $H_1(C;\ZZ)$.
Since $j_{\#}$ is trivial on vanishing cycles, the lemma follows.
\qed
\msk
By the above lemma, $j_{\#}$ induces a surjective homomorphism
$\ov{j}_{\#}\cl\pi_1(V)\to\pi_1(\cJ)$. We will finish the proof by showing that
$\ov{j}_{\#}$ is trivial. This is a consequence of the following easy result.

\proclaim (2.3.3) Lemma.
Let $D\ss\CC$ be a disc centered at $0$, and let $\psi\cl D\hra U$ be a (holomorphic)
embedding such that:
\item{\rm 1.}  
$\psi(D^0)\ss V$, where $D^0:=D\sm\{0\}$,
\item{\rm 2.}
$\psi(0)$ is a smooth point of $\D$ (or $\L$), and $D$ intersects $\D$ (respectively $\L$)
transversely.
\n
Then, after shrinking $D$, we can assume there exists a lift $\wt{\psi}\cl D\to \cJ$ of
$\psi$. 

\pf
Let $q\cl\cC\to D$ be the family of curves in $|\cO_X(2)|$ parametrized by $D$. By Item~(2)
the analytic surface $\cC$ is smooth. Shrinking $D$ we can assume there exists a
line-bundle $\cL$ on $\cC$ such that for each $t\in D$, the couple
$(q^{-1}t,\cL|_{q^{-1}t})$ satisfies the hypotheses of Lemma~(2.1.2). By the modular
property of $\cJ$ the couple $(\cC,\cL)$ induces a morphism $\wt{\psi}\cl D\to\cJ$
lifting $\psi$. 
\qed
\msk
To finish the proof that $\ov{j}_{\#}$ is trivial, let
$R\ss U|$ be a straight line (notice that $|\cO_X(2)|\sm U$ has codimension $3$) transverse
to $\D$ and $\L$. The map induced by inclusion
$$\pi_1(R\sm (\D\cup\L))\to \pi_1(V)$$
is a surjection by Bertini's Theorem~[GM, p.151]. The fundamental group on the left is
generated by $\pi_1(D^0_{x_i})$, where $R\cap (\D\cup\L)=\{x_1,\ldots,x_k\}$, and
$D^0_{x_i}\ss(R\sm (\D\cup\L))$ is a small punctured disc centered at $x_i$. By
Lemma~(2.3.3) the map $\ov{j}_{\#}$ is trivial on each  $\pi_1(D^0_{x_i})$. This shows
$\ov{j}_{\#}$ is trivial, and concludes the proof that $\cJ$ is simply-connecetd. 

\vfill
\eject
\n
{\bf 3. Proof that $b_2(\cMt_4)\ge 24$.}
\bsk
\n
Let $\phi\cl\cM_4\to\PP^N$ be the  determinant map~(1.1). Composing Donaldson's
map $\mu\cl H_2(X;\QQ)\to H^2(\phi(\cM_4);\QQ)$ (see~[FM,Mo])
with pull-back by $\phi$, we get 
$$\phi^{*}\mu\cl H_2(X;\QQ)\to H^2(\cM_4;\QQ).$$
This map is an injection, because by~[FM,Mo]
$$\int_{\cM_4}\wedge^{10}(\phi^{*}\mu(\a))=
{10!\over 5!2^5}\left(\int_X\wedge^2\a\right)^5.\leqno(3.1)$$
Consider the boundary divisor $B\ss\cM_4$ (see~(1.1)). By Lemma~(1.9.1) $B$ contains
$\PP^1$'s contracted by $\phi$.  Hence, if $h\in
H^2(\cM_4)$ is the class of a hyperplane, $\QQ h$ and
$\phi^{*}\mu(H_2(X;\QQ))$ span a $23$-dimensional subspace of
$H^2(\cM_4)$. Now we pull back by $\wt{\pi}\cl\cMt_4\to\cM_4$. 

\proclaim (3.2) Claim.
Keep notation as above. Then
$$\dim\left(\QQ\pit^{*}h+\pit^{*}\phi^{*}\mu(H_2(X;\QQ))\right)=23.$$ 

\pf
The subspace   $\wt{\pi}^{*}\phi^{*}\mu(H_2(X;\QQ))$ has dimension $22$ by~(3.1). If
$\wt{\PP}^1\ss\cMt_4$ is the proper transform of a $\PP^1\ss B$ contracted by $\phi$,
then  
$$\langle \pit^{*}h,\wt{\PP}^1\rangle\not=0,\qquad 
\langle \wt{\pi}^{*}\phi^{*}\mu(H_2(X;\QQ)),\wt{\PP}^1\rangle=0.$$
Thus  $\QQ\pit^{*}h$ is not contained in
$\pit^{*}\phi^{*}\mu(H_2(X;\QQ))$.
\qed
\msk
Finally let's show that $\QQ c_1(\Sit)$ 
and $\QQ\pit^{*}h+\pit^{*}\phi^{*}\mu(H_2(X;\QQ))$ span a
$24$-dimensional subspace of $H^2(\cMt_4;\QQ)$. Let $\Sit_{Z,W}\ss\cMt_4$ be the
fiber of $\pit$ over $[I_Z\op I_W]$, where $Z\not= W$. By~[O4, (2.3.1)]
$\Sit_{Z,W}\cong\PP^1$, and by~[O4,(2.3.2)]
$$\langle c_1(\Sit),\Sit_{Z,W}\rangle =-2,\qquad 
\langle \QQ\pit^{*}h\op\pit^{*}\phi^{*}\mu(H_2(X;\QQ)),\Sit_{Z,W}\rangle=0.$$
Thus $\QQ c_1(\Sit)$ 
is not contained in $\QQ\pit^{*}h+\pit^{*}\phi^{*}\mu(H_2(X;\QQ))$. By Claim~(3.2) we
conclude that
$$\dim\left(\QQ c_1(\Sit)+\QQ\pit^{*}h+\pit^{*}\phi^{*}\mu(H_2(X;\QQ))\right)=24.$$  
\bsk
\centerline{\bf References.}
\msk
\item {[B]} A.~Beauville. {\it Vari\'et\'es K\"ahl\'eriennes dont la
premi\`ere classe de Chern est nulle}, J.~Differential Geom.~18 (1983),
755-782.
\item {[FM]} R.~Friedman-J.~Morgan. {\it Smooth four-manifolds and
complex surfaces}, Ergeb.~Math.~Grenzgeb. (3. Folge) 27, Springer (1994).
\item{[GH]} L.~G\"ottsche, D.~Huybrechts. {\it Hodge numbers of moduli
spaces of stable bundles on $K3$ surfaces}, preprint, MPI/94-80. 
\item {[GM]} M.~Goresky, R.~MacPherson. {\it Stratified Morse Theory}, 
Ergeb.~Math.~Grenzgeb. (3. Folge) 14, Springer (1988). 
\item {[H]} D.~Huybrechts. {\it Compact Hyperk\"ahler manifolds: basic
results}, alg-geom/9705025. 
\item {[La]} R.~Lazarsfeld. {\it Brill-Noether-Petri without degenerations}, J.~of
Differential Geom.~23 (1986), 299-307.  
\item {[LP1]} J.~Le Potier. {\it Fibr\'e d\'eterminant et courbes de saut sur les
surfaces alg\'ebriques}, Complex projective geometry, London
Math.~Soc.~Lecture Note Series 179, Cambridge University Press (1992).
\item {[LP2]} J.~Le Potier. {\it Syst\`emes coh\'erents et structures de niveau},
Ast\'erisque 214, Soc.~Math.~de France (1993).
\item {[Li1]} J.~Li. {\it Algebraic geometric interpretation of Donaldson's
polynomial invariants of algebraic surfaces}, J.~Diff.~Geom.~37 (1993),
417-466.
\item {[Li2]} J.~Li. {\it The first two Betti numbers of the moduli spaces of
vector bundles on surfaces}, preprint (1995). 
\item {[M1]} S.~Mukai. {\it Symplectic structure of the moduli space of
sheaves on an abelian or $K3$ surface}, Invent.~math.~77 (1984), 101-116.
\item {[M2]} S.~Mukai. {\it On the moduli space of bundles on $K3$ surfaces,
I}, in Vector bundles on algebraic varieties, T.I.F.R., Oxford Univ.~Press
(1987), 341-413.
\item {[O1]}  K.~G.~O'Grady. {\it The weight-two Hodge  structure of moduli
spaces of sheaves on a K3 surface}, J.~of Algebraic Geom.~6 (1997),
599-644. 
\item {[O2]}  K.~G.~O'Grady. {\it Moduli of vector-bundles on surfaces},
Algebraic Geometry Santa Cruz 1995,
Proc.~Symp. Pure Math.~vol.~62, Amer.~Math.~Soc.~(1997), 101-126.
\item {[O3]} K.~O'Grady. {\it Moduli of vector bundles on projective
surfaces: some basic results}, Invent.~math.~123 (1996), 141-207.
\item{[O4]} K.~O'Grady.  {\it Desingularized moduli spaces of sheaves on a
$K3$, I}, preprint Dip.to di Matematica ``G.~Castelnuovo'' 98/21 (1998).
\item{[R]} M.~Raynaud. {\it Sections des fibr\'es vectoriels sur une courbe},
Bull.~Soc.~Math.~Fr.~110 (1982), 103-125.  
\item{[S]} C.~Simpson. {\it Moduli of representations of the fundamental
group of a smooth projective variety I}, Publ.~Math.~Inst.~Hautes
\'Etudes Sci.~79 (1994), 47-129.   
\bsk
\bsk
\n 
Universit\'a di Roma ``La Sapienza'' 
\msk
\n 
Dipartimento di Matematica ``G.~Castelnuovo''
\msk
\n
00185 Roma
\msk
\n 
ITALIA
\bsk
\n
e-mail: ogrady@mat.uniroma1.it

\end